\documentclass[12]{article}
\usepackage[centertags]{amsmath}
\usepackage{amsfonts}
\usepackage{amssymb}
\usepackage{amsthm}
\usepackage{newlfont}
\usepackage[ansinew]{inputenc}
\usepackage[english,francais]{babel}
\usepackage{amstext}

\theoremstyle{plain}
\newtheorem{Theorem}{Theorem}[section]
\newtheorem{proposition}{Proposition}[section]
\newtheorem{Corollary}{Corollary}[section]
\newtheorem{Definition}{Definition}[section]
\newtheorem{Lemma}{Lemma}[section]
\newtheorem{Remark}{Remark}[section]
\newtheorem{Conjecture}{Conjecture}[section]

\begin{document}

\title{The $A$-Module Structure Induced by a Drinfeld $A$-Module of Rank 2 over a Finite Field }
\author{{\normalsize {MOHAMED AHMED Mohamed Saadbouh}} \\
%EndAName
{\small {Institut de Mathématiques de Luminy, CNRS-UPR 9016} }\\
{\small {Case 907, 163 Avenue de Luminy, 13288 Marseille, Cedex
9. France.}}}
\date{}
\maketitle

\selectlanguage{francais}
\begin{abstract}
Soit $\Phi $ un $\mathbf{F}_{q}[T]$-module de Drinfeld de rang
$2$, sur un
corps fini $L$, une extension de degré $n$ d'un corps fini à $q$ éléments $%
\mathbf{F}_{q}$. \\Soit $P_{\Phi }(X)=$ $X^{2}-cX+\mu P^{m}$( $c$ un élément de $%
\mathbf{F}_{q}[T]$, $\mu $ un élément non nul de
$\mathbf{F}_{q}$, $m$ le degré de
l'extension $L$ sur le corps $\mathbf{F}_{q}[T]/P$ et $P$ est la $\mathbf{F}%
_{q}[T]$-caractéristique de $L$) le polyn%
ôme caractéristique, de Frobenius $F$ de $L$. On s'inté%
ressera à la structure de $\mathbf{F}_{q}[T]$-module fini
$L^{\Phi }$ et on
prouvera notre résultat principal qui est l'analogue du théorème de Deuring pour les courbes elliptiques:\ \ \ soit $M=\frac{\mathbf{F}%
_{q}[T]}{I_{1}}\oplus \frac{\mathbf{F}_{q}[T]}{I_{2}}$, où $I_{1}=(i_{1})$,$%
\ \ I_{2}=(i_{2})$\ ( $i_{1}$, $i_{2}$ deux polynômes de $\mathbf{F}_{q}[T]$%
) et tel que : $i_{2}\mid (c-2)$. \\Alors il existe un $\mathbf{F}_{q}[T]$%
--module de Drinfeld $\Phi $ sur $L$ de rang $2$, ordinaire, tel que : $%
L^{\Phi }$ $\simeq M$. On fini par une statistique concernant la
cyclicité de la structure de $A$-module $L^{\Phi }$.
\end{abstract}

\selectlanguage{english}

\begin{abstract}
Let $\Phi $ be  a Drinfeld $\mathbf{F}_{q}[T]$-module of rank
$2$, over a finite field $L$, a finite extension of $n$ degrees
for a finite field with $q$ elements $%
\mathbf{F}_{q}$. Let $P_{\Phi }(X)=$ $X^{2}-cX+\mu P^{m}$ ($c$ an element of $%
\mathbf{F}_{q}[T]$ and $\mu $ a no null element of
$\mathbf{F}_{q}$, $m$  the degree of the
extension $L$ over the field  $\mathbf{F}_{q}[T]/P$, $P$ is a $\mathbf{F}%
_{q}[T]$-characteristic of $L$ and $d$ the degree of the
polynomial $P$) the
 characteristic polynomial, of the Frobenius $F$ of $L$. We will interested to the structure of finite
$\mathbf{F}_{q}[T]$-module  $L^{\Phi }$ deduct by $\Phi $ over
$L$ and will proof our main result, the analogue of Deuring theorem for the elliptic curves  :\ \ \ Let $M=\frac{\mathbf{F}%
_{q}[T]}{I_{1}}\oplus \frac{\mathbf{F}_{q}[T]}{I_{2}}$, where $I_{1}=(i_{1})$,$%
\ \ I_{2}=(i_{2})$\ ( $i_{1}$, $i_{2}$ two polynomials of $\mathbf{F}_{q}[T]$%
) and such that : $i_{2}\mid (c-2)$. Then there exists an ordinary Drinfeld  $\mathbf{F}_{q}[T]$%
-module $\Phi $ over $L$ of rank $2$, such that : $%
L^{\Phi }$ $\simeq M$. We finish by a statistic about the
cyclicity of such structure $L^{\Phi }$, and we prove that is
cyclic only for the trivial extensions of $\mathbf{F}_{q}$.
\end{abstract}

\section{Introduction}
Let $E$ be an elliptic curve over finite field $\mathbf{F}_{q}$,
we know,
see [12], [15], [16] and [19], that the endomorphism ring of $E$, End$_{%
\mathbf{F}_{q}}E$, is an  $\mathbb{Z}-$order in  a division
algebra, this algebra is  : $\mathbb{Q}$ and in this case End$_{\mathbf{F}_{q}}E=\mathbb{Z}$%
, or a quadratic complex field and in this case : End$_{\mathbf{F}_{q}}E=%
\mathbb{Z}+c$ $O_{K}$ where $c$ is an element of $\mathbb{Z}$ and
$O_{K}$ is the maximal $\mathbb{Z}-$order in this quadratic
complex field, or is a Quarternion Algebra over $\mathbb{Q}$ and
in this case  End$_{F_{q}}E$ is a maximal order in this
Quarternion  Algebra. We put  $E(\mathbf{F}_{q}) $ the abelian
group of \textbf{\ }$\mathbf{F}_{q}$-rational points of $E$. The
cardinal of this abelian  group is equal to $N=q+1-c$, and by
Hasse-Weil $\mid c\mid \leq 2.\sqrt{q}$. The structure of this
group in the ordinary case is  :

\begin{equation*}
E(\mathbf{F}_{q})\simeq \mathbb{Z}/A\oplus \mathbb{Z}/B,\text{ if
}(c,q)=1,\text{ \ }B\mid A,\text{ \ }B\mid (c-2)\text{ \ \ and
}A.B=N\qquad (1)
\end{equation*}%

conversely, for every abelian group  $\mathbb{Z}/A\oplus \mathbb{%
Z}/B,$ with $B\mid A$, $B\mid (c-2)$ \ \ and $A.B=N,$ there is an
elliptic curve $E$ such that : $\ $ $E(\mathbf{F}_{q})\simeq
\mathbb{Z}/A\oplus \mathbb{Z}/B,$ we note that in the
supersingular case, this structure is also known, see [15]. \\Over
this structure, S.\ Vladut, in [20], has effected a statistic
about the report of elliptic curves for
 which the $E(\mathbf{F}_{q})$ is cyclic over the number of classes of $\mathbf{F}_{q}$-endomorphisms of elliptic curves over
finite field $\mathbf{F}_{q}$, this report will depend on $q$ and
will be noted  $c(q)$, and we have :
\begin{equation*}
c(q)=\frac{\text{\#\{}E\text{, }E(\mathbf{F}_{q})\text{ cyclic\}}}{\#\{E\}}%
,
\end{equation*}%
where $\#\{E\}$ is the  number of classes of
$\mathbf{F}_{q}$-endomorphisms of elliptic curve over a finite
field $\mathbf{F}_{q}$, and we know, see  [20], that :

$c(q)=1$ if and only if  $q=2^{l}$ where $l\neq 2$ is a prime  or
equal at $1$ and  one of the following conditions is satisfied:

\begin{enumerate}
\item $q-1$ is prime, $q\neq 4$ ( the case $q=2$ is included, thus we consider $1$
as prime ),

\item $q-1=l_{1}l_{2}$ with primes  $l_{1}$ and $l_{2}$ is not a
\textquotedblright\ small \textquotedblright\ divisors of $q-1$;

\item $\ q-1=l_{1}l_{2}l_{3}$ \ with primes  $l_{1},l_{2}$ and $%
l_{3} $ are not a  \textquotedblright\ small \textquotedblright\
divisors of $q-1$.
\end{enumerate}
the case $l_{1}=l_{2}$ is not exclude.

In general case, the number $c(q)$ is given in [20], by :

Let $\varepsilon > 0$ we have:
\begin{equation*}
c(q)=\prod\limits_{l}(1-\frac{1}{l(l^{2}-1)})+O(q^{-1/2+\varepsilon
}),
\end{equation*}
where the product is taken over all prime divisors of  $q-1$. Our
goal here is to give an analog of the above mentioned results in
the case of Drinfeld Modules of rank 2. We recall quickly what is
it : let $K$ a no empty global field of characteristic $p$ ( that
means a rational functions field of one indeterminate over a
finite field ) with a constant field the finite field
$\mathbf{F}_{q}$ with $p^{s}$ elements. We fix one place of  $K$,
noted $\infty $ \ and we call  $A$ \ the ring of regular elements
away from the place $\infty $. Let $L$ be a commutative field of
characteristic $p$, and let $\gamma :A\rightarrow L,$ be a
$A$-ring homomorphism, the kernel of this  homomorphism is noted $P$ and $m$\ =$[L,$ $%
A/P]$ is the extension degrees of  $L$ over $A/P$.

We note $L\{\tau \}$ the Ore's  polynomial ring, that means the
polynomial ring of  $\tau $, $\tau $ is  the Frobenius of
$\mathbf{F}_{q}$, with the usual  addition and the product is
given by the computation rule : for every   $\lambda $ of
$L\mathit{,}$ $\tau \lambda =\lambda ^{q}\tau $. We say a
Drinfeld $A$-module $\Phi $ for a non trivial homomorphism of
ring, from $A$ to \ $L\{\tau \}$ which is different of $\gamma $.
this homomorphism $\Phi $, once defined, gives a  $A$-module
structure over the  $A$-field $L$, noted \ $L^{\Phi }$, where the
name of a Drinfeld $A$-module for a homomorphism $\Phi $. This
structure $A$-module is depending on $\Phi $ and especially on
this rank.

Let $\chi $ be the characteristic of Euler-Poincare ( it is a
ideal from $A$ ), so we can speak about the ideal  $\chi (L^{\Phi
})$, will be noted by $\chi _{\Phi }$, which is by Definition a
divisor for $A,$ corresponding for the elliptic curves to a
number of points of the variety over their basic fields.

We will work, in this paper, in the special case $K=$ $\mathbf{F}_{q}(T)$, $\ A=%
\mathbf{F}_{q}[T]$. Let $P_{\Phi }(X):$ the characteristic polynomial of the $A$%
-module $\Phi $, it is also a characteristic polynomial of
Frobenius $F$ of $L$.\ We can prove that this polynomial can be
given by : $P_{\Phi }(X)=$ $X^{2}-cX+\mu P^{m}$, such that $\mu
\in \mathbf{F}_{q}^{\ast } $, $c$ $\in A$ and $\ \deg c\leq
m.\frac{d}{2}$, the Hass-Weil's analogue in this case.

We  will  interest  to a Drinfeld $A$-module structure $L^{\Phi }$
in the case of rank 2, and we will prove that for an ordinary
Drinfeld $\mathbf{F}_{q}[T]$-module, this structure is always the sum of two cyclic and finite $\mathbf{F}_{q}[T]$%
-modules : $\frac{A}{I_{1}}\oplus \frac{A}{I_{2}}$ where  $%
I_{1}=(i_{1})$ and $I_{2}=(i_{2})$ such that $i_{1}$ and $i_{2}$
is two ideals of $A,$ which verifies   $i_{2}\mid i_{1}$. We will
show that  $\chi _{\Phi }=I_{1}I_{2}=(P_{\Phi }(1))$, and if we
put \ $i=$paced$(i_{1},i_{2}),$ then : $i^{2}\mid P_{\Phi }(1)$. \
We will give now some appears of our results proved in this paper
:
\begin{proposition}
With the above notations, we have :\newline $L^{\Phi }\simeq
\frac{A}{I_{1}}\oplus \frac{A}{I_{2}}$. And if we have an
ordinary module $\Phi $, then \ : $\ i_{2}\mid (c-2)$.
\end{proposition}

We note  by  End$_{L}\Phi $ the endomorphism ring of a Drinfeld $A$-module $%
\Phi $, we have:

\begin{proposition}
Let $\Phi $ \ be  an ordinary Drinfeld $A$-module of rank 2 and
let
$\rho $ a prime ideal of $A$ different from $P$ the $A$-characteristic of $%
L$, such that $\rho ^{2}\mid P_{\Phi }(1)$ and $\rho \mid (c-2)$.
Then $\Phi (\rho )\subset L^{\Phi }$ if and only if the $A$-order
$\ O(\Delta /\rho ^{2})\subset $ \ End$_{L}\Phi $.
\end{proposition}

Finally, we come  to our main result, which is a complete analog
of (1), the Deuring-Waterhous theorem  for the elliptic curves :

\begin{Theorem}
Let $M=\frac{A}{I_{1}}\oplus \frac{A}{I_{2}}$, where
$I_{1}=(i_{1})$ ,$\ \ I_{2}=(i_{2})$\ and such that : $\
i_{2}\mid i_{1}, $ $i_{2}\mid (c-2)$. Then there exists an
ordinary Drinfeld $\ \mathit{A}$-module  $\Phi $ over $L$ of rank
$2$, such that: $L^{\Phi }\simeq M$.
\end{Theorem}

Lastly, we will make  a statistic about the ordinary Drinfeld
$A$-modules such that the $A$-modules $L^{\Phi }$ are cyclic, we
note by $C(d,m,q)$ the proportion of the number (of isomorphisms
of) ordinary Drinfeld \ $A$-modules, of rank 2 such that the
A-modules structures $L^{\Phi }$ are cyclic, this means : if we
note by \\ $ \#\{\Phi$, isomorphism, ordinary \} the  number of
classes of $L$-isomorphisms of \\an ordinary Drinfeld Modules of
rank 2, we have :

\begin{equation*}
C(d,m,q)=\frac{\#\{\Phi \text{, }L^{\Phi }\text{ cyclic}\}}{\#\{\Phi \text{%
, isomorphism, ordinary}\}}\text{,}
\end{equation*}
and we note by  $C_{0}(d,m,q)$ the proportion of the number
( of isogeny classes of) ordinary Drinfeld $A$-modules, of rank  $2$ such that the $A$%
-modules $L^{\Phi }$ are cyclic, otherwise, if we note by $%
\#\{\Phi $, isogeny, ordinary$\}$ the number of isogeny classes,
of ordinary Drinfeld modules of rank $2$, we have :
\begin{equation*}
C_{0}(d,m,q)=\frac{\#\{\text{isogeny Classes of }\Phi \text{, }L^{\Phi }%
\text{ cyclic}\}}{\#\{\Phi \text{, isogeny, ordinary}\}}\text{.}
\end{equation*}

Of course, these numbers are depending on $q$ and also on $d,m$.\\
One of our important results is :

\begin{proposition}
$C(d,m,q)=C_{0}(d,m,q)=1$if and only if $m=d=1$.
\end{proposition}

This means that, to have a Cyclic Drinfeld  $A$-modules we must
have a trivial extension $L$, we give also some values for  $C(d,m,q)$ and $%
C_{0}(d,m,q)$ corresponding to some given values of $d$ and $m$,
for example :

\begin{proposition}
We put $d=2$ and $m=1$. Let $H(O(D))$ the Hurwitz's number of
classes for an order $O$ which the imaginary determinant is $D$ :
\begin{equation*}
C_{0}(2,1,q)=\frac{q(q-1)-5}{q(q-1)-2}\text{,}
\end{equation*}
\end{proposition}

\begin{equation*}
C(2,1,q)=\frac{q^{3}-q^{2}-q+1-[\frac{q-1}{2}\sum\limits_{P_{\Phi
}}\sum\limits_{i_{2,}i_{2}^{2}\shortmid 4-4\mu P}H(O(\frac{4-4\mu P}{%
i_{2}^{2}}))+(q-1)\sum\limits_{P_{\Phi
}}\sum\limits_{i_{2},i_{2}^{2}\shortmid c^{2}-4\mu P}H(O(\frac{c^{2}-4\mu P}{%
i_{2}^{2}}))]}{q^{3}-q^{2}-q+1}\text{.}
\end{equation*}%
And we  let think, in conjecture form, that for a big $q$ the
values of $C(d,m,q)$ and $C_{0}(d,m,q)$ will tend to 1.

\bigskip
\bigskip

\section{Drinfeld Modules }

Let $E$ be an extension of $\mathbf{F}_{q}$, and let $\tau $ Frobenius of $%
\mathbf{F}_{q}$. We put $E\{\tau \}$ the polynomial ring in $\tau
$ with the usual  addition and the multiplication defined by:$\ $

\begin{equation*}
\ \forall e\in E,\text{ }\tau e=e^{q}\tau \text{.}
\end{equation*}
\begin{Definition}

Let  $R$ be  the $E$-linearly polynomials set with the
coefficient in
$E$, that means that these elements are on the following form :%
\begin{equation*}
Q(x)=\sum_{K>0}l_{k}x^{q^{K}}\text{,}
\end{equation*}%
where $l_{k}\in E$ for every $k>0$, and only a finite number of
$l_{k}$ is not null. The ring  $R$ is a ring by addition and the
polynomial composition.
\end{Definition}

\begin{Lemma}
$E\{\tau \}$ and $R$ are isomorphic rings.
\end{Lemma}

If we put  \ $A=\mathbf{F}_{q}[T],$ $\ f(\tau
)=\sum_{i=0}^{v}a_{i}\tau ^{i}\in E\mathit{\{\tau \}}$ \textit{and
\ }$Df\mathit{:=a}_{0}=f^{\prime }(\tau ).$

It is clear that the application :
\begin{equation*}
E\mathit{\{\tau \}\mapsto }E
\end{equation*}

\begin{equation*}
f\mapsto Df,
\end{equation*}%
is a  morphism of $\mathbf{F}_{q}$-algebras.
\begin{Definition}
 An $A$-fields $E$ is a field $\mathit{E}$
equipped with a fix morphism \\ $ \gamma :A \longrightarrow E$.
The prime ideal $P=$ Ker$\gamma $ is called the characteristic of
$E$. We say $E$ has generic characteristic if and only if
$P=(0)$; otherwise (i.e $P \neq (0)$) we said $P=$ is finite and
$E$ has finite characteristic.
\end{Definition}

We then have the following fundamental definition :

\begin{Definition}
Let $E$ an $A$-field and $\Phi :A\mapsto E\mathit{\{\tau \}}$ be
homomorphism of algebra. Then $\Phi $ is a Drinfeld $A$-module
over $E$ if and only if  :

\begin{enumerate}
\item $D\circ \Phi =\gamma ;$

\item for some $a$ $\in A,\Phi _{a}\neq \gamma (a)\tau ^{0}.$
\end{enumerate}
\end{Definition}

As was proved by Drinfeld in [6], such modules always exist.
\begin{Remark}

\begin{enumerate}
\item The above normalization is analogous  to the normalization used in complex
multiplication of elliptic curves. The last condition is
obviously a non-triviality condition.
\item By $\Phi $, every  extension $E_{0}$ of $E$ became an  $A$-module
by :
\begin{eqnarray*}
E_{0}\times A &\rightarrow &E_{0}, \\
\ (k,a) &\rightarrow &k.a:=\Phi _{a}(k)\ .
\end{eqnarray*}%
We will note this  $A$-module by $E_{0}^{\Phi }.$
\end{enumerate}
\end{Remark}

Let $\overline{E}$ be a fix algebraic closure of $E$ and $\Phi $ a
Drinfeld module over $E$ and $I$ an ideal of $A$. As $A$ is a
Dedekind domain, one know that $I$ may be generated by ( at most ) two elements $%
\{i_{1}$ , $i_{2}\}\subset I.$

Since $E\mathit{\{\tau \}}$ has a right division algorithm, there
exists a right greatest common divisor in $E\mathit{\{\tau \}}$.
It is the monic generator of the left ideal of $E\mathit{\{\tau
\}}$ generated by : $\Phi _{i_{1}}$ et $\Phi _{i_{2}}.$

\begin{Definition}
We set  $\Phi _{I}$ to be the monic generator of the left ideal
of $E\mathit{\{\tau \}}$ generated by   $\Phi _{i_{1}}$ and
$\Phi_{i_{2}}.$
\end{Definition}

\begin{Definition}
Let $E_{0}$ an extension of  $E$ and $I$ an ideal of $A$. We
define by  $: $ $\ \Phi \lbrack I](E)$ the finite subgroup of $\Phi \lbrack \overline{%
E_{0}}]$ given by the roots of $\Phi _{I}$ in $\overline{E}$.
\end{Definition} If  $a\in A$, then we set $\Phi \lbrack a]:=\Phi
\lbrack (a)].$ We can see it as :\\
 $\Phi \lbrack a]=\{$ set of roots of $\Phi (a)$ in $\overline{K}%
\},$ and $\Phi _{I}=\cap _{a\in I}\Phi \lbrack a]$. Then :

\begin{equation*}
\Phi _{a}(\overline{E}):=\Phi \lbrack a](\overline{E})=\{x\in \overline{E}%
,\Phi _{a}(x)=0\}.
\end{equation*}

And for every ideal  $Q\subset A,$
\begin{equation*}
\Phi _{Q}(\overline{E}):=\Phi _{Q}(\overline{E})=\cap _{a\in Q}\Phi _{a}(%
\overline{E}).
\end{equation*}

\begin{Remark}
The groups \ : $\Phi \lbrack I](E)$ and $\Phi \lbrack
I](\overline{E})$ are clearly stable under  $\{\Phi _{a}\}_{a\in
A}$.\ \ \
\end{Remark}

\begin{Definition}
Let $\Phi $ be a Drinfeld $A$-module over an  $A$-field $E$. We
say that  $ \Phi $ is supersingular, if and only if, the
$A$-module constituted by a $P$-division points $\Phi
_{P}(\overline{E})$ is trivial.
\end{Definition}

%\end{Definition}

\subsection{The Height and Rank of a Drinfeld Module $\Phi $}

Let $\Phi $ be a Drinfeld $A$-module over the $A$-field $E$. We
note by deg$_{\tau }$ the degree in indeterminate $\tau $.

\begin{Definition}
An element of $E\{\tau \}$ is called separable, if this constant
coefficient is not null. It called purely inseparable if it is on
the form  $\lambda \tau ^{n}$, $n>1$ and $\lambda \in E$, $\lambda
\neq 0$ .
\end{Definition}

Let $H$ be a global field of characteristic $p>0$, and let
$\infty $ one place (a Prime ideal ) of $H$, we will note by
$H_{\infty }$ the completude of $H$ at the place  $\infty $. We
define the degree of function over $A$ by :

\begin{Definition}
Let $a$ $\in A,$ deg $a=\dim _{F_{q}}\frac{A}{aA}$ if $a$ $\neq 0$
and $\deg 0=$ $-\infty $.
\end{Definition}

We extend  deg to $\ K$ by putting deg $x=$ deg $a-$deg $b$ if $0\neq x$ $%
= $ $\frac{a}{b}$ $\in K$.  \\If $A=\mathbf{F}_{q}[T]$, then the
degree function is the usual polynomial  degree. Let $Q$ be a no
null ideal  of  $A$, we define the ideal degrees of $%
Q$, noted deg $Q$, by :
\begin{equation*}
\text{deg}Q=\dim _{F_{q}}\frac{A}{Q}.
\end{equation*}

\begin{Lemma}
there exists a rational number $r$ such that  :
\begin{equation*}
\deg _{\tau }\Phi _{a}=r\deg a\text{.}
\end{equation*}
\end{Lemma}
\begin{proof}
It is easy to see that $\Phi $ is an injection, otherwise since
$K\{\tau \}$ is an integre ring, Ker $\Phi $ is a prime ideal non
null, so maximal in $A$ and by consequence  Im $\Phi $ is a
field, so $\Phi =\gamma $. Since $-$deg$_{\tau }$ define  a no
trivial valuation over Frac$(\Phi (A))$ ( the fractions field of
$\Phi (A)$ ) which is isomorphic to $K$, so $-$deg$_{\tau }$ and
$-$deg are  equivalent valuations over $K$. There is rational
number $r>0 $, such that :
\begin{equation*}
r\ \text{deg}=\text{deg}_{\tau }.
\end{equation*}
\end{proof}
\begin{Corollary}
Let $\Phi :A\mapsto E\{\tau \}$ be a Drinfeld $A$-module, so
$\Phi $ is injective.
\end{Corollary}

\begin{proposition}
The number  $r$ is a positive integer.
\end{proposition}

\begin{Definition}
The number $r$ is called the rank of the Drinfeld $A$-module
$\Phi.$
\end{Definition}

For example if $A=\mathbf{F}_{q}[T],$ a  Drinfeld $A$-module of
rank $r$ is on the  form  : $\ \Phi (T)=a_{1}+a_{2}\tau
+...+a_{r}\tau ^{r}$, where  $\ a_{i}\in E$, $1\leq i\leq r-1$
and $a_{r}\in E^{\ast }$.

In this case  char $E\mathbf{=}$ $P\neq (0)$ we can define the
notion of height of a Drinfeld module $\Phi$.

For this,  we put  $\ v_{P}:K\mapsto
\mathbb{Z},$ an associate normalized valuation at $P$, this means, if  $a$ $%
\in K$ has a root over  $P$ of order $t$, we have $v_{P}(a)=t$.

For every $a$ $\in A,$ let $w(a)$ the most small integer $t>0$,
where $\tau ^{t}$ occurs  at $\Phi _{a}$ with a no null
coefficient.

\begin{Lemma}
There exists a rational number $h$ such that :%
\begin{equation*}
w(a)=hv_{P}(a)\deg P\text{.}
\end{equation*}
\end{Lemma}

\begin{proposition}
The number $h$ is a positive integer.
\end{proposition}

\begin{Definition}
the integer $h$ is called the height of $\Phi .$
\end{Definition}

For example if $A=\mathbf{F}_{q}[T],$ a Drinfeld $A$-module of height $%
h $ of rank  $r$ is on the form : $\Phi (T)=a_{0}+a_{h}\tau
^{h}+...+a_{r}\tau ^{r},$ where $\ a_{i}\in E$, $0\leq i\leq r-1$
and $a_{r}\in E^{\ast }$.

\begin{Definition}
Let $\Phi $ and $\Psi $ two Drinfeld $A$-modules over an
$A$-field $E$ and $p$ an isogeny over  $E$ from $\Phi $ to $\Psi
.$
\end{Definition}

\begin{enumerate}
\item We say that $p$ is separable if and only if  $p(\tau )$ is
separable.
\item We say that  $p$ is purely no separable if and only if $p(\tau
)=\tau ^{j}$ for one $j>0$.
\end{enumerate}

\subsection{Norm of Isogeny}

\begin{Definition}
Let $F$ an integer over a ring $A,$ with fractions field $K$. we
note by $N_{K/K(F)}$ the determinant of the $K$-linearly
application of multiplication by  $F$ to $K(F)$ ( it is the usual norm if the extension $%
K(F)/K$ is  separable.
\end{Definition}

We can see that there is a morphism  $N_{K/K(F)}:I_{\overline{A}
}\rightarrow I_{A\text{ }}$ from the fractional ideals groups of  $\overline{%
A}$ to functionary ideals group of $A,$ by this morphism we have:

\begin{proposition}
The norm of isogeny is a principal ideal.
\end{proposition}

\begin{proposition}
Let $M_{_{fin}}(A)$ the category of primes ideals of $A$ and let
$D(A)$ the  monoïde of ideals of  $A$. There exists an unique
function  :
\begin{equation*}
\chi : M_{_{fin}}(A)\mapsto D(A),
\end{equation*}
multiplicative over the exact sequence and such that $\chi (0)$
$=1$ and $\chi (A/\wp $ $)=\wp $ for every prime ideal $\wp $ of
$A$.
\end{proposition}

\begin{Definition}
The function $\chi $ is called the Euler-Poincare characteristic.
\end{Definition}

We can see $\chi (E^{\Phi })$ and we note it by  $\chi _{\Phi }.$

\begin{proposition}
The ideals $\chi _{\Phi }$ and \ $P^{m}$ are principals (in $A$),
and more clearly  $\chi _{\Phi }=(P_{\Phi }(1))$ and
$P^{m}=P_{\Phi }(0)$.
\end{proposition}

\begin{enumerate}
\item We know that the norm of isogeny is a principal ideal, indeed  $%
N(F)=$ $P_{\Phi }(0)$ and $N(1-F)=$ $(P_{\Phi }(1))$ since $F$ and
$1-F$ are a  $K$-isogenys.

\item We can call  $\chi _{\Phi }$ the divisor of  $E$-points,
this divisor is  analogue at the number of $E$-points for elliptic
curves.

\item $\chi _{\Phi }$ is the annulator of  $A$-module $E^{\Phi }$. We can deduct that  : $E^{\Phi }\subset (\frac{A}{\chi _{\Phi }})^{r}.$

\item The structure of $A$-module $E^{\Phi }$ is stable by the Frobenius endomorphism $F$.
\end{enumerate}
\begin{Corollary}
If there are a Drinfeld  $A$-module  $\Phi$, over a field  $E,$ of
characteristic $P$ and of degree  $m$ over $A/P$, then the ideal
$P^{m}$ is a principal ideal.
\end{Corollary}

\begin{Remark}
The above Corollary shows that there exists  a  restriction of
the existence of Drinfeld $A$-modules.
\end{Remark}

\section{Drinfeld Modules over Finite Fields}

We substitute the Field $E$, by $L$ a finite extension of degree
$n$ of the finite field $\mathbf{F}_{q}$. Let $\tau :x\mapsto
x^{q}$ be the  Frobenius of $\mathbf{F}_{q}$, so the Frobenius
$F$ of $L$ is \ $F=\tau ^{n}$ and $\mathbf{F}_{q}[F]$ is the
center of $L\{\tau \}$. We put  $m=[L:A/P]$ and $d=$ deg $P$, then
$n=m.d$. The function $-$deg define a  valuation over $K$, the
field of fractions of $A$. We put  $r_{1}=[K(F):K]$ and
$r_{2}^{2}$ the degree of the left field  End$_{L}\Phi \otimes
_{A}K$ over this center  $K(F)$.

So a Drinfeld $A$-module $\Phi $ over $L$ give a  structure of
$A$-module over the additive finite group $L$, this  structure
will be noted  $L^{\Phi }$. Let $\gamma $ the application of $A$
to $L$ which an element $a$ for $A$ associate the constant of
$\Phi _{a}$, then it is easy to see that  $\gamma $ is a ring
homomorphism, and that $\Phi $ and
 $\gamma $ are equal over  $A^{\ast }=\mathbf{F}_{q}^{\ast }$ the set of reversible elements
 of $A$.

\begin{Definition}
Let $\Phi $ be a Drinfeld $A$-module over a finite field $L$. We note by  $%
M_{\Phi }(X)$ the unitary minimal polynomial of $F$ over $K$ .
\end{Definition}

\begin{proposition}
With the above notations  : $M_{\Phi }(X)$ is an element of
$A[X]$, equal to $P_{\Phi }^{\frac{1}{r_{2}}}.$
\end{proposition}

\begin{Corollary}
For two Drinfeld $A$-modules  $\Phi $ and $\Psi $, of rank $r$
over a finite field $L$, then the following are equivalent :

\begin{enumerate}
\item $\Phi $ and $\Psi $ are isogenous,

\item $M_{\Phi }(X)=M_{\Psi }(X),$

\item $P_{\Phi }=P_{\Psi }.$
\end{enumerate}
\end{Corollary}

\begin{proposition}
Let $L$ be a finite extension of degree  $n$ over a finite field
$\mathbf{F}_{q}$ and $F$ the  Frobenius of $L$. Then $L(\tau )$ is
a central division algebra over $\mathbf{F}_{q}(F)$ of dimension
$n^{2}.$
\end{proposition}
\begin{Definition}
Every  $u$ $\in L\{\tau \}$ can be writing on this  form $u=\tau
^{h}u^{\prime }$ ( since $L$ is a perfect field) where $u^{\prime
}\in L\{\tau \}$ separable. The integer $h$ is called the height
of  $u$ and will be noted by ht $u.$
\end{Definition}
In the finite field case, we can see the height of a Drinfeld
$A$-module\ $\Phi $ over finite field $L$, the integer $H_{\Phi
}$, as  been :
\begin{equation*}
H_{\Phi }=\frac{1}{r}\inf \{\text{ht }\Phi _{a},0\neq a\in P\}.
\end{equation*}

\begin{Remark}
It is easy to see that  $H_{\Phi }$ is invariant under isogeny
and that
\begin{equation*}
1\leq H_{\Phi }\leq r.
\end{equation*}
\end{Remark}

\begin{proposition}
Let  $\Phi $ be a Drinfeld  $A$-module of rank  $r$ over a finite
field  $L$, the following assertions are equivalent :

\begin{enumerate}
\item There exists a finite extension  $L^{\prime }$ of $L,$ such that
the  endomorphism ring End$_{L^{\prime }}\Phi \otimes _{A}K,$ has
dimension $r^{2}$ over  $K$.

\item Some power of the  Frobenius $F$ of $L$ lies in $A.$

\item $\Phi $ is  supersingular.

\item The field  $K(F)$ has only one prime above $P$.
\end{enumerate}
\end{proposition}
\begin{proposition}
Let $\Phi $ be a Drinfeld  A-module of rank  $r$ and let  $Q$ be
an ideal from $A$ prime with $P,$ then :%
\begin{equation*}
\Phi _{Q}(\overline{L})=(\frac{A}{Q})^{r}\text{.}
\end{equation*}
\end{proposition}
\begin{Corollary}
Then we can deduct that  : $\Phi
_{P}(\overline{L})=(\frac{A}{P})^{r-H_{\Phi }} $.
\end{Corollary}

We can deduct from above mentioned Proposition the following
important result, which characterize the supersingularity :

\begin{proposition}
\ The Drinfeld $A$-module $\Phi $ is supersingular ( $\Phi
_{P}(\overline{L})=0$ ), if and only if, $r=H_{\Phi }$.
\end{proposition}
\begin{Definition}
We say that the field  $L$ is so big if all endomorphism rings
defined over  $\overline{L}$
are also defined over $L$, i.e : End$_{\overline{L%
}}\Phi =$ End$_{L}\Phi .$
\end{Definition}

Two Drinfeld modules $\Phi $ and $\Psi $ are isomorphic, if and
only if, there exists an $a$ $\in $ $L$ such that : $a^{-1}\Phi
_{a}=\Psi _{a}a.$
\begin{Lemma}
Let $\Phi $ be a Drinfeld $ A$-module of rank $r$, over a finite
field $L$, of characteristic $P$. The characteristic polynomial
 of Frobenius endomorphism $F$ is :%
\begin{equation*}
P_{\Phi }(X)=X^{r}+c_{1}X^{r-1}+...+c_{r-1}X+\mu P^{m},c_{1},...c_{r-1}\in A%
\text{ et }\mu \in \mathbf{F}_{q}^{\ast }\text{.}
\end{equation*}
\end{Lemma}

\begin{Remark}
The fact that constant of the polynomial $P$ is  $\mu P^{m}$ comes
from the fact that $P_{\Phi }(0)=P^{m}$ in $A.$
\end{Remark}

The following Proposition is an analogue of the Riemann's
hypothesis for elliptic curves :

\begin{proposition}
Let $\Phi $ be a Drinfeld $A$-module of rank $r$ over finite
field $L$ which is a finite extension of degree $n$ of \
$\mathbf{F}_{q}$. Then deg $(w)=\frac{n}{r}$\, for every root $w$
of characteristic polynomial $P_{\Phi }(X).$
\end{proposition}

The following result is the Hasse-Weil's analogue for the elliptic
curves :
\begin{Corollary}
Let $P_{\Phi }(X)=X^{r}+c_{1}X^{r-1}+...+c_{r}X$ $+\mu P^{m}$ be
the  characteristic polynomial of a Drinfeld Module $\Phi$, of
rank  $r,$ over a finite field $L$. Then:

\begin{equation*}
\forall 1\leq i\leq r-1,\text{deg }c_{i}\leq \frac{i}{r}m\deg P.
\end{equation*}
\end{Corollary}
\begin{proof}
The  proof can be deducted immediately by the above Proposition.
\end{proof}
\section{Drinfeld Modules of rank 2}

In all next of this paper, $\Phi $ will be considered a Drinfeld
$A$-module of rank $2$, And $A=\mathbf{F}_{q}[T]$ for proof and
more details see [1], [12] and [6].

\subsection{Structure of A-module $L^{\Phi }$}

Let $\Phi $ be a Drinfeld $A$-module of rank $2,$ \ over a finite
field  $L$ and let $P$ this characteristic polynomial. About the
$A$-module structure $L^{\Phi }$, we have the following result :

\begin{proposition}
The Drinfeld $A$-module $\Phi $ give a finite $A$-module structure
$L^{\Phi }$,  which is on the form  $\frac{A}{I_{1}}\oplus \frac{A}{I_{2}}$ where $%
I_{1}$ and $I_{2}$ are two ideals of $A,$ such that: $\chi _{\Phi
}=I_{1}I_{2}.$
\end{proposition}

\begin{proof}
Since the $A$-module  $L^{\Phi }$ is a sub$A$-module of $\Phi
(\chi _{\Phi })$ $\simeq $ $\frac{A}{\chi _{\Phi }}\oplus
\frac{A}{\chi _{\Phi }}$, then there are $I_{1}$ and $I_{2}$ in
$A$ such that  : $L^{\Phi }\simeq \frac{A}{I_{1}}\oplus
\frac{A}{I_{2}}$ and since the Euler-Poincare's Characteristic is
multiplicative over the exact sequence  we will have $\chi _{\Phi
}=I_{1}I_{2}.$
\end{proof}

We put  $I_{1}=(i_{1})$ and $\ \ I_{2}=(i_{2})$ ( $i_{1}$ and
$i_{2}$ two unitary polynomials in $A$).

Let \ $i=$ pgcd $(i_{1},i_{2})$, it is clear by the Chinese lemma,
that the non cyclicity of the $A$-module $L^{\Phi }$, needs that
\ $I_{1}$ and $I_{2}$ are not a prime between them, that means
that $i\neq 1,$ and since the relation $\chi _{\Phi
}=I_{1}I_{2}$, we will have  : $\ i^{2}\mid P_{\Phi }(1)$ $\
(\chi _{\Phi }=(P_{\Phi }(1))$).

In all the next of this paper, the condition above,  will be
considered verified, and more precisely we suppose that $I_{2}\mid
I_{1}$ $\ ($i.e $:i_{2}\mid i_{1})$ otherwise $L^{\Phi }$ is a
cyclic  $A$-module and can be writing on this form  $ A/\chi
_{\Phi }.$

\begin{proposition}
If $L^{\Phi }\simeq \frac{A}{I_{1}}\oplus \frac{A}{I_{2}}$, then
$ i_{2}\mid c-2$.
\end{proposition}

\begin{proof}
We know that the  $A$-module structure $L^{\Phi }$ is  stable by
the endomorphisme  Frobenius $F$ of $L$. We choose a basis for $ A
/\chi _{\Phi },$ for which  the $A$-module
$L^{\Phi }$ will be generated by $(i_{1},0)$ and $(0,i_{2})$.\\
Let $M_{F}$ $\in \mathbf{M}_{2}(A/\chi _{\Phi })$ the  matrix of
the endomorphisme  Frobenius $F$ in this basis. Then
$M_{F}=~\left(
\begin{array}{ll}
a & b \\
a_{1} & b_{1}
\end{array}
\right)$, where $a,b,a_{1},b_{1}\in A/\chi _{\Phi }.$\\
 Although since :
Tr $M_{F}=a+b_{1}=c$ and $M_{F}(i_{1},0)=(i_{1},0)$ and $
M_{F}(0,i_{2})=(0,i_{2})$, we will have  $a.i_{1}\simeq
i_{1}(\mod\chi _{\Phi })$ and then  $a-1$ is divisible by
$i_{1}$, of same for  $ b_{1}.i_{2}\simeq i_{2}(\mod \chi _{\Phi
})$, that means that  $ b_{1}-1$ is divisible by $i_{2}$ and
then: $c-2=a-1+b_{1}-1$ is divisible by $i_{2}$ ( since we have
always  $i_{2}\mid i_{1}$).
\end{proof}

Let \ $\rho $ be a prime ideal from  $A$, different from the  $A$ -characteristic $P$%
, we define the finite $A$-module $\Phi (\rho )$ as been the $A$-module $%
(A/\rho )^{2}.$

The discriminant of the $A$-order: $A+g.O_{K(F)}$  is  $\Delta .g^{2}$, where $%
\Delta $ is the discriminant of the characteristic polynomial
$P_{\Phi }(X)=X^{2}-cX+\mu P^{m}$. So each order is defined by
this discriminant and will be noted by $O($ disc$)$. It is clear,
by the Propositions 4.1  that the inclusion $\Phi (\rho )\subset
$ $L^{\Phi }$
implies that $%
\rho ^{2}\mid P_{\Phi }(1)$ and $\rho \mid c-2$. And if we note by End$%
_{L}\Phi $ the endomorphism ring of the Drinfled $A$-module $\Phi
$, we have :
\begin{proposition}
Let $\Phi $ be an ordinary Drinfeld  $A$-module of rank 2, and let
$\rho $ an ideal from  $A$ different from the $A$-characteristic $P$ of $L,$ such that $%
\rho ^{2}\mid P_{\Phi }(1)$ and $\rho \mid c-2$. Then$\Phi (\rho
)\subset L^{\Phi }$, if and only if, the $A$-order $\ O(\Delta
/\rho ^{2})\subset $ End$_{L}\Phi $.
\end{proposition}

To prove this Proposition we need the following lemma :

\begin{Lemma}
$\Phi (\rho )\subset L^{\Phi }$ is equivalent to
$\frac{F-1}{\varrho }\in $ End$_{L}\Phi $.
\end{Lemma}

\begin{proof}
Since $L^{\Phi }=$ Ker$(F-1)$ and $\Phi (\rho )$ $=$ Ker$(\rho )$
( We confuse  by commodity the ideal $\rho $ with this generator
in $A$) and we know by [3], Proposition 4.7.9, that for two
isogenys, let by example $F-1$ and $\rho $, we have
 Ker$(F-1)\subset $ \
Ker$(\rho )$, if and only if, there exists an element $g\in $ End$%
_{L}\Phi $ such that $F-1=g.\rho $ and then $\Phi (\rho )\subset
L^{\Phi },$ if and only if, $\frac{F-1}{\varrho }=g\in
$End$_{L}\Phi .$
\end{proof}

We prove now the Proposition 4.3 :

\begin{proof}
Let $N(\frac{F-1}{\rho })$ the norm of the isogeny
$\frac{F-1}{\rho }$, which is a principal ideal generated by
$\frac{P_{\Phi }(1)}{\rho ^{2}}$, and the trace $($Tr$)$ of this
isogeny is $\frac{c-2}{\rho }$ then we can calculate the
discriminant of the  $A$-module $\ A[\frac{ F-1}{\rho }]$ by:

disc$A([\frac{F-1}{\rho }])=Tr(\frac{F-1}{\rho
})^{2}-4N(\frac{F-1}{\varrho } )=\frac{c^{2}-4\mu P^{m}}{\rho
^{2}}=\Delta /\rho ^{2}\Rightarrow$

$$O(\Delta /\rho ^{2})\subset  \ End_{L}\Phi. $$

We suppose now that : $O(\Delta /\rho ^{2})\subset $ End$_{L}\Phi
$ and we prove that $\Phi (\rho )\subset L^{\Phi }$. The Order
corresponding of the discriminant  $ \Delta /\rho ^{2}$ is
$A[\frac{F-1}{\rho }]$ this means that: $\frac{ F-1}{\varrho }\in
$ End$_{L}\Phi $ and so, by lemma 4.1 : $\Phi (\rho )\subset
L^{\Phi }$ .
\end{proof}

\begin{Corollary}
If $O(\Delta /\rho ^{2})\subset $ \ End$_{L}\Phi $, then $L^{\Phi
}$ is not cyclic.
\end{Corollary}

\begin{proof}
We know that $\Phi (\rho )$ is not cyclic  (since it is a
$A$-module of rank  $2$), and then the necessary and sufficient
conditions need for non cyclicity of $A$-module \ $L^{\Phi }$ are
equivalent to the necessary and  sufficient conditions to have
$\Phi (\rho )\subset L^{\Phi }$.
\end{proof}

We can so prove the following important Theorem :

\begin{Theorem}
Let $M=\frac{A}{I_{1}}\oplus \frac{A}{I_{2}}$, $I_{1}=(i_{1})$
And $\ \ I_{2}=(i_{2})$\, such that : $i_{2}\mid i_{1},$\
$i_{2}\mid (c-2)$. Then there exists an ordinary Drinfeld
$\mathit{A}$-module $\Phi $ over $L$ of rank $2$, such that:
\begin{equation*}
L^{\Phi }\simeq M.
\end{equation*}
\end{Theorem}
\begin{proof}
In fact, if we consider the Drinfeld $A$-module $\Phi$, for which
the characteristic of Euler-Poincare is giving by $\chi _{\Phi
}=I_{1}.I_{2}$ and this endomorphism ring is $O(\Delta
/i_{2}^{2})$ where $\Delta $ is always the discriminant of the
characteristic polynomial of the Frobenius $F$. We remind   that\
$\Phi (\rho )\subset L^{\Phi }$ for every  $\rho $ an ideal $A$,
different from $P$ and verify  $\rho ^{2}\mid P_{\Phi }(1)$ and
$\rho \mid (c-2),$ if and only if, the $A$-order $\ O(\Delta /\rho
^{2})\subset $ \ End$_{L}\Phi $. Let now  $\rho =i_{2}$. Since by
construction the $A$-order $\ O(\Delta /i_{2}^{2})\subset $
End$_{L}\Phi $ we have that $\Phi (i_{2})\simeq
(A/i_{2})^{2}\subset L^{\Phi }$. We know that $L^{\Phi }$ is
included or equal to $\Phi (\chi _{\Phi })$ $\simeq $
$\frac{A}{\chi _{\Phi }} \oplus \frac{A}{\chi _{\Phi }}$, we have
so  : $L^{\Phi }=\frac{ A}{I_{1}}\oplus \frac{A}{I_{2}}$.
\end{proof}

The Theorem above can be proved by using the following conjecture:

\begin{Conjecture}
Let $M\in \mathbf{M}_{2}(A/\chi _{\Phi })$ , \
$\overline{P}=P(\mod \chi
_{\Phi })$. \\We suppose : $($det $M)=\overline{P}^{m}$, Tr $(M)=c$ and $c$ $%
\nmid $ $P$. There exists a ordinary Drinfeld  $A$-module over a
finite field $L$ of rank $2$, for which the Frobenius matrix
associated, is  $M_{F}$, and such that  : $\ $
\begin{equation*}
M_{F}=M\in \mathbf{M}_{2}(A/\chi _{\Phi })\text{.}
\end{equation*}
\end{Conjecture}

We put the following matrix :
\begin{equation*}
M_{F}=\left(
\begin{array}{ll}
c-1 & i_{1} \\
i_{2} & -1%
\end{array}
\right) \in \mathbf{M}_{2}(A/\chi _{\Phi })\text{.}
\end{equation*}

We can see that the three  conditions of the  conjecture are
realized then there exists an ordinary Drinfeld $A$-modules $\Phi
$ over $L$ of rank $2$, such that : $L^{\Phi }\simeq M$.

\subsection{Deuring Theorem }

The following Theorem, proved by Max-Deuring in [15] and [19] is
used for the proof of the analogue of our principal result, in
elliptic curves case :

\begin{Theorem}
Let $E_{0}$ be an elliptic curve over a finite field of
characteristic  $p,$ with a no trivial endomorphism $F_{0}$. Then
there exists an elliptic curve $E$ over a field of numbers and an
endomorphisme $F$ from $E$ such   $E_{0}$ is isomorphic to
$\overline{E}$ and $F_{0}$ corresponding to $F$ under this
isomorphism.
\end{Theorem}

From the previous Theorem, we can deduct the following Theorem :

\begin{Theorem} let $N$ $\in \mathbb{N}$, $M=\left(
\begin{array}{ll}
a & b \\
a_{1} & b_{1}%
\end{array}
\right) \in \mathbf{M}_{2}( \mathbb{Z}/N\mathbb{Z})$ \ and
$\mathbf{F}_{q}$ a finite field with $q$ elements, we suppose:

\begin{enumerate}
\item $($det $M)=q(\mod N);$

\item $\mid a+b_{1}\mid \leq 2.\sqrt{q}.$

There exists a Frobenius endomorphism $F$ which verifies : $
F^{2}-cF+q=0$ (mod N),  such that $c=a+b_{1}$ and this matrix
$M_{F}\in \mathbf{M}_{2}( \mathbb{Z}/N\mathbb{Z})$ is exactly $M$.
\end{enumerate}
\end{Theorem}

This Theorem is used to prove  the following Theorem :

\begin{Theorem}
Let $M=\left(
\begin{array}{ll}
c-1 & -A \\
B & 1%
\end{array}
\right) \in \mathbf{M}_{2}(\mathbb{Z}/N\mathbb{Z} )$ and
$\mathbf{F}_{q}$ is a finite field with $q$ elements, such that :
$\mid c\mid $ $\leq 2.\sqrt{q}$, $B\mid A $, $B\mid c-2$
$A.B=N=q+1-c,$ we suppose : \ $(c,q)=1$. Then there exists an
ordinary elliptic curve $E$ over $\mathbf{F}_{q}$, such that :
\end{Theorem}

\begin{equation*}
E(\mathbf{F}_{q})\simeq \mathbb{Z}/A\oplus \mathbb{Z}/B\text{.}
\end{equation*}

\section{Cyclicity Statistics for the $A$-module L$^{\Phi }$}

In this section, we make a statistic about the Drinfeld Modules
$\Phi$ of rank 2, whose the structures \ $L^{\Phi }$ are cyclic, \
for this, we define  $C(d,m,q)$ as been the ration of the number
of (isomorphism classes of)  Drinfeld modules of rank 2 with
cyclic structure $L^{\Phi }$ to the number of $L$-isomorphisms
classes of ordinary Drinfeld modules of rank  $2$, noted by
$\#\{\Phi $, isomorphism, ordinary$\}$:

\begin{equation*}
C(d,m,q)=\frac{\#\{\Phi ,L^{\Phi }\text{ cyclic}\}}{\#\{\Phi
\text{, isomorphism, ordinary}\}}\text{,}
\end{equation*}

As same, we define   $N(d,m,q)$ as been the ration of the number
of (isogeny classes of)  Drinfeld modules of rank 2 with not
cyclic structure $L^{\Phi }$ to the number of $L$-isomorphisms
isogeny of ordinary Drinfeld modules of rank  $2$, noted by
$\#\{\Phi $, isogeny, ordinary$\}$:

\begin{equation*}
N(d,m,q)=\frac{\#\{\Phi ,L^{\Phi }\text{ non cyclic}\}}{\#\{\Phi
\text{ , isogeny, ordinary}\}}\text{.}
\end{equation*}

We remark that : $\ 0\leq C(d,m,q)$, $N(d,m,q)\leq 1$. Since the
no cyclicity of the structure $L^{\Phi }$ needs the fact that  $\
i^{2}\mid P_{\Phi }(1)$ and $i_{2}\mid (c-2)$, it is natural to
introduce $i$ ( so  $i_{2}$ ) in the calculus of  $C(d,m,q)$ and
$N(d,m,q)$.

We fix the characteristic polynomial $P_{\Phi }$, this means that
we fix the isogeny classes of $\Phi $, and we define :

\begin{Definition}
We note by $n(P_{\Phi },i_{2})=\#\{\Phi :L^{\Phi
}=\frac{A}{(i_{1})}\oplus \frac{A}{(i_{2})}\}$.
\end{Definition}

\begin{Remark}
The number $n(P_{\Phi },i_{2})$ is equal to the number of
isomorphisms classes of $\Phi$ whole the $A$-module $L^{\Phi }\simeq \frac{A}{%
(i_{1})}\oplus \frac{A}{(i_{2})}$, in one isogeny classes, from
where is coming the  correspondence between $\Phi $ and $i_{2}.$
\end{Remark}

For $n(P_{\Phi }, i_{2})$ we have by the Theorem 3.1 :

\begin{Lemma} Let $P_{\Phi }(X)=X^{2}-cX+\mu P^{m}$ be the characteristic
polynomial of an ordinary Drinfeld $A$-module $\Phi$ of rank 2,
and let $i_{2}$ be an unitary polynomial of $A$. Then if$\ $
 $i_{2}\mid c-2$ we have : $n(P_{\Phi },i_{2})\geq 1$, else
$n(P_{\Phi },i_{2})=0$.
\end{Lemma}
We can deduct  :

\begin{Corollary}
with the above notations :
\begin{equation*}
\#\{\Phi ,L^{\Phi }\text{ non cyclic}\}=\sum\limits_{P_{\Phi
}}\sum\limits_{i_{2},i_{2}^{2}\mid P_{\Phi }(1)}n(P_{\Phi
},i_{2}).\#\{i_{2},i_{2}^{2}\mid P_{\Phi }(1)\text{ and }i_{2}\mid
(c-2)\} \text{,}
\end{equation*}

\begin{equation*}
\#\{\Phi ,L^{\Phi }\text{ cyclic}\}=\sum\limits_{P_{\Phi
}}\sum\limits_{i_{2},i_{2}^{2}\mid P_{\Phi }(1)}n(P_{\Phi
},i_{2}).\#\{i_{2},i_{2}^{2}\nmid P_{\Phi }(1)\text{ and \
}i_{2}\mid (c-2)\} \text{;}
\end{equation*}

And if we note by $n_{0}(P_{\Phi },i_{2})=$ $\#\{$isogeny classes
of $ \Phi :L^{\Phi }=\frac{A}{(i_{1})}\oplus \frac{A}{(i_{2})}\}$
, we have :
\begin{equation*}
n_{0}(P_{\Phi },i_{2})=1\text{.}
\end{equation*}
\end{Corollary}
We note now by  $\#\{\Phi ,$ isogeny, ordinary$\}$ the number of
isogeny classes, for an ordinary module  $\Phi $, then we define :
\begin{equation*}
N_{0}(d,m,q)=\frac{\#\{\text{isogeny classes of }\Phi \text{, }L^{\Phi }%
\text{ not cyclic}\}}{\#\{\Phi ,\text{isogeny,
ordinary}\}}\text{,}
\end{equation*}

the same for
\begin{equation*}
C_{0}(d,m,q)=\frac{\#\{\text{isogeny Classes of }\Phi \text{, }L^{\Phi }%
\text{ cyclic}\}}{\#\{\Phi ,\text{isogeny, ordinary}\}}\text{,}
\end{equation*}

We can so announce the following lemma :

\begin{Lemma}
With the notations above, we have :
\begin{eqnarray*}
N_{0}(d,m,q) &=&\frac{\text{\#}\{i_{2},i_{2}^{2}\mid P_{\Phi }(1)\text{ and}%
i_{2}\mid (c-2)\}}{\#\{\Phi \text{, isogeny, ordinary}\}}\text{,} \\
N(d,m,q) &=&\frac{\sum\limits_{P_{\Phi
}}\sum\limits_{i_{2},i_{2}^{2}\mid P_{\Phi }(1)}n(\Phi
,i_{2}).\text{\#}\{i_{2},i_{2}^{2}\mid P_{\Phi }(1)\text{ and
}i_{2}\mid (c-2)\}}{\#\{\Phi \text{, isomorphism,
ordinary}\}}\text{,}
\end{eqnarray*}

\begin{equation*}
C_{0}(d,m,q)=\frac{\text{\#}\{i_{2},i_{2}^{2}\nmid P_{\Phi
}(1)\text{ et } i_{2}\mid (c-2)\}}{\#\{\Phi \text{, isogeny,
ordinary}\}}\text{,}
\end{equation*}%
\begin{equation*}
C(d,m,q)=\frac{\sum\limits_{P_{\Phi
}}\sum\limits_{i_{2},i_{2}^{2}\nmid P_{\Phi }(1)}n(\Phi
,i_{2}).\text{\#}\{i_{2},i_{2}^{2}\nmid P_{\Phi }(1) \text{ and
}i_{2}\mid (c-2)\}}{\#\{\Phi \text{, isomorphism, ordinary}\}}
\text{,}
\end{equation*}

and $N(d,m,q)+C(d,m,q)=1$, $N_{0}(d,m,q)+C_{0}(d,m,q)=1$.

\end{Lemma}

The calculus of $\#\{\Phi ,$ isogeny, ordinary$\}$, for an
ordinary $A$-module $\Phi $, has been calculated in  [3] and [4],
as been :

\begin{proposition}
Let $L=F_{q^{n}}$ and $P$ the $A$-characteristic of $L$.

We put $m=[L:A/P]$ and $d=$deg $P$ :

\begin{enumerate}
\item $m$ is odd and $d$ is odd :%
\begin{equation*}
\#\{\Phi ,\text{ isogeny, ordinary}\}=(q-1)(q^{[\frac{m}{2}d]+1}-q^{[\frac{%
m-2}{2}d]+1}+1).
\end{equation*}
\item $m.d$ even:
\begin{equation*}
\#\{\Phi ,\text{ isogeny, ordinary}\}=(q-1)(\frac{(q-1)}{2}q^{\frac{m}{ 2}%
d}-q^{\frac{m-2}{2}d}+1)\text{.}
\end{equation*}
\end{enumerate}
\end{proposition}

As for the number $L$-isomorphisms classes, we will need the
following result, for the proof and more details see [9] :

\begin{proposition}
Let $L$ be a finite extension of degree $\ n$ over
$\mathbf{F}_{q}$, then the number of $L$-isomorphisms classes of
a Drinfeld $A$-module of rank  2 over $L$ is $(q-1)q^{n}$ if $n$
is odd and $q^{n+1}-q^{n}+q^{2}-q$ else.
\end{proposition}

And to calculate the number of  $L$-isomorphisms classes for an
ordinary Drinfeld  modules, we will need to calculate the number
of $L$-isomorphisms classes for an supersingular Drinfeld modules
and subtract it from the global number of $L$-isomorphisms
classes, for this, we have by [10] :

\begin{proposition}
Let $L$ be a finite extension of $\ n$ degrees over
$\mathbf{F}_{q}$, then the number of   $L$-isomorphisms classes
of an supersingular Drinfeld $A$-module of rank 2, over $L$ is
$(q^{n_{2}}-1)$, where $n_{2}=$ pgcd($2,n $).
\end{proposition}

The calculus of $C(d,m,q)$ will be calculated in function of the
values of $d$ and $m$ which are two major values to determinate
$c$ because  deg $c\leq \frac{m.d}{2}$.

And to calculate the number of $L$-isomorphisms classes existing
in each isogeny classes, we need the following Definition for
more information, see [13]:

\begin{Definition}
Let $L$ be a finite extension of degree $\ n$ over
$\mathbf{F}_{q},$ we define $W(F)$ as been :
\begin{equation*}
W(F)=\sum\limits_{\Phi ,\text{ }F=\text{Frobenius(}\Phi
\text{)}}Weigh(\Phi )
\end{equation*}%
where :%
\begin{equation*}
\text{Weigh}(\Phi )=\frac{q-1}{\#\text{Aut}_{L}\Phi }\text{.}
\end{equation*}
\end{Definition}

$W(F)$ is the sum of weights( noted Weigh$(\Phi )$ ) of number of
$L$-isomorphisms classes existing in each isogeny classes of the
module $\Phi $ which the  Frobenius is $F$.

And to calculate $\#$Aut$_{L}\Phi $ we have the following lemma :

\begin{Lemma}
Let $\Phi $ be an ordinary Drinfeld $A$-module of rank $2$, over a
finite field  $L=F_{q^{n}} $, then : $\#$Aut$_{L}\Phi $ = $q-1$.
\end{Lemma}

By the previous lemma, we can see that Weight $(\Phi )=\frac{ q-1}{\#%
\text{Aut}_{L}\Phi }$ $=1$, that means :

\begin{Corollary}
In the case of ordinary Drinfeld modules of rank  2, $W(F)$ is
the number of  $L$-isomorphisms classes existing in each isogeny
classes.
\end{Corollary}
\begin{Definition}
Let $D$ be an imaginary discriminant and let $l$ a polynomial for
which  the square is a divisor of $D$ and let
$h(\frac{D}{l^{2}})$ the number of classes of the order for which
the discriminant is $\frac{D}{l^{2}}$. We define the number of
classes of Hurwitz for an imaginary discriminant  $D$, noted
$H(D) $ by :
\begin{equation*}
H(D)=\sum\limits_{l}\sum\limits_{l^{2}\shortmid
D}h(\frac{D}{l^{2}})\text{.}
\end{equation*}
\end{Definition}

For more definitions and information about the numbers classess
of Hurwitz, see [13] and [20].

\begin{Lemma}
If $\alpha $ is an integral  element over $A$, for which $O
=A[\alpha ]$ is an $A$-order, then disc$(A[\alpha ])$ is equal to
the discriminant of the minimal polynomial of $\alpha .$
\end{Lemma}

What is  interesting for us is the calculus of the disc$(A[F])$
and since disc$(A[F])=$ disc$(P_{\Phi })$. To calculate the number
of classes  $W(F)$, we have the following result, for proof see
[13].

\begin{proposition}
Let $L$ be \ a finite extension of degree $n$ of a field  $F_{q}$
and $F$ the Frobenius of $L$, then:

\begin{equation*}
W(F)=H(\text{disc}(A[F]))\text{.}
\end{equation*}
\end{proposition}
It remains for us to calculate  $n(\Phi ,i_{2}):$

\begin{Lemma}
Let $P_{\Phi }$ be the characteristic polynomial of an ordinary
Drinfeld $A$-module of rank $2$, over a finite field $L$ such that $L^{\Phi }= \frac{A}{%
(i_{1})}\oplus \frac{A}{(i_{2})}$, and let  $\Delta $ the
discriminant of the characteristic polynomial of the Frobenius
$F$, then :
\begin{equation*}
n(P_{\Phi },i_{2})=H(O(\Delta /i_{2}^{2}))\text{.}
\end{equation*}
\end{Lemma}

\begin{proof}

We know that to have  $L^{\Phi }=\frac{A}{(i_{1})}\oplus
\frac{A}{(i_{2})}$, we have  certainly :

$\Phi (i_{2})\simeq (A/i_{2})^{2}\subset L^{\Phi }$, that is
equivalent to say, by the Proposition 4.3, that the $A$-order $\
O(\Delta /i_{2}^{2})\subset $ End$_{L}\Phi $ where $\Delta $ is
always the  discriminant of the characteristic polynomial of $F$,
$P_{\Phi }$, and that :
\begin{equation*}
n(P_{\Phi },i_{2})=H(O(\Delta /i_{2}^{2}))\text{.}
\end{equation*}
\end{proof}

We can calculate some values of $C(d,m,q)$ for some $d$ and $%
m$ :

\section{The case : $d=m=1$}

In this case $L=A/P=\mathbf{F}_{q}$, the $A$-module $L^{\Phi }$
$=A/P$, so it is  cyclic, that means that  $C(1,1,q)=1$.

We can find this result by more explicit way  :

\begin{proposition}
Let $L=F_{q^{n}}$ and $P$ the $A$-characteristic of $L$,
$m=[L,A/P]$ \\and $ d= $ deg$P$. We  suppose $m=d=1$. then :
\begin{equation*}
C(1,1,q)=C_{0}(1,1,q)=1\text{.}
\end{equation*}
\end{proposition}

\begin{proof}
$P_{\Phi }(1)=1-c+\mu P^{m}=1-c+\mu P$, \ : $\ i_{2}^{2}\mid
P_{\Phi }(1)$ $\Rightarrow $ $i_{2}$ is a no null constant, then
an element of $\mathbf{F}_{q}^{\ast }$, so  $(i_{2})=A$ and
$\frac{A}{I_{1}} \oplus \frac{A}{I_{2}}=\frac{A}{I_{1}}$, then
$L^{\Phi }$ is  cyclic, that means that $C(1,1,q)=1\Rightarrow
N(1,1,q)=0$.

To calculate \#\{$i_{2},i_{2}^{2}\mid P_{\Phi }(1)$ and $i_{2}\mid
(c-2)$\}, we must have  $i_{2}$ an element of
$\mathbf{F}_{q}^{\ast }$ and deg $ i_{2}>0\Rightarrow $

\#\{$i_{2},i_{2}^{2}\mid P_{\Phi }(1)$ and $i_{2}\mid (c-2)$\}$=0$
and then :

\begin{equation*}
N_{0}(1,1,q)=\frac{\#\{i_{2},i_{2}^{2}\mid P_{\Phi }(1)\text{ and
}i_{2}\mid (c-2)\}}{\#\{\Phi \text{, isogeny, ordinary}\}}=0
\end{equation*}%
$\Rightarrow C_{0}(1,1,q)=1$.
\end{proof}

By more precise way, we can announce :

\begin{Theorem}
Let $L=F_{q^{n}}$ and $P$ the $A$-characteristic of $L$,
$m=[L,A/P]$ and $ d=$ deg$P$. Then:
\begin{equation*}
C_{0}(d,m,q)=C(d,m,q)=1\Leftrightarrow m=d=1\text{.}
\end{equation*}
\end{Theorem}

\begin{proof}

We have just seen that  $C_{0}(1,1,q)=$ $C(1,1,q)=1$.

Conversely and by the absurd  : $m.d>1$ (this means $%
m.d\geq 2)$, we take for example $m=1$ and $d=2$ .

To have $C_{0}(d,m,q)=$ $C(d,m,q)=1$, we must have

$\#\{i_{2},i_{2}^{2}\mid P_{\Phi }(1)$ and $i_{2}\mid (c-2)\}=0$,
what is not true, since if $c=aT+b$, where $a\in
\mathbf{F}_{q}^{\ast }$ and $b\in
\mathbf{F}_{q}$, it is sufficient to have an  unitary $i_{2}$ and such that : $%
i_{2}\mid (c-2)$, for this, we take : $i_{2}=a^{-1}(c-2)$ this
stay compatible with the fact that $a^{-2}(c-2)^{2}\mid 1-c+\mu
P$, since there are many solutions for the equation in $i_{2}$,
i.e in $a$ and $b:$ \

\begin{equation*}
a^{-2}(c-2)^{2}\mid 1-c+\mu P\Rightarrow a^{-2}(aT+b-2)^{2}\mid
1-aT-b+\mu (T^{2}+pT+p_{0}),
\end{equation*}

which implies the equations : $2a^{-1}\mu (b-2)=\mu p_{1}-a$ and
$\mu \lbrack a^{-1}(b-2)]^{2}=1-b+\mu p_{0}\Rightarrow 2^{-1}[\mu
p_{1}-a]=1-b+\mu p_{0}$ this gives a values of $a$ for each value
of  $b$, from where the many possibilities  of  $i_{2}$, for
example $i_{2}=$ $T-(\mu (2p_{0}-p_{1}))^{-1}$, so it is more
easy in the case which $m.d>2$ to find an $i_{2},$ such that :
$i_{2}^{2}\mid P_{\Phi }(1)$ and $i_{2}\mid (c-2)$.
\end{proof}
\section{The case: $m=1$ and $d=2$}

In this case $n=m.d=2$, and $n_{2}=2\Rightarrow \#\{\Phi $,
isomorphism, ordinary$\}=$ $q^{3}-q-(q^{2}-1)=q^{3}-q^{2}-q+1.$

\begin{proposition}
Let $L=F_{q^{n}}$ and $P$ the $A$-characteristic of $L$,
$m=[L,A/P]$ and $\ d=$ deg $P$. We suppose  $m=1$ and $d=2$. Then:

\begin{equation*}
C_{0}(2,1,q)=\frac{q(q-1)-5}{q(q-1)-2}\text{,}
\end{equation*}

\begin{equation*}
C(2,1,q)=\frac{q^{3}-q^{2}-q+1-[\frac{q-1}{2}\sum\limits_{P_{\Phi
}}\sum\limits_{i_{2,}i_{2}^{2}\shortmid 4-4\mu P}H(O(\frac{4-4\mu P}{%
i_{2}^{2}}))+(q-1)\sum\limits_{P_{\Phi
}}\sum\limits_{i_{2},i_{2}^{2}\shortmid c^{2}-4\mu P}H(O(\frac{c^{2}-4\mu P}{%
i_{2}^{2}}))]}{q^{3}-q^{2}-q+1}\text{.}
\end{equation*}
\end{proposition}
\begin{proof}
 We start by calculating :
\begin{equation*}
\frac{\#\{i_{2},i_{2}^{2}\mid P_{\Phi }(1)\}}{\#\{\Phi \text{;
isogeny} \}}.
\end{equation*}
\ For this, we distinguish between two cases, the case where $c=2$
and the case where  $c \neq 2$ .

Then for $c=2$ : $\ i_{2}^{2}\mid P_{\Phi }(1)\Rightarrow
i_{2}^{2}\mid \mu P^{m}-1$ this implies that if we put
$i_{2}=T+j_{2}$, $j_{2}\in \mathbf{F}_{q}$ and
$P(T)=T^{2}+p_{1}T+p_{0}$ where $p_{1}$, $p_{0}\in
\mathbf{F}_{q}$, are irreducibles, we will have $p_{1}=2j_{2}$ and
$\mu p_{0}-1=\mu j_{2}^{2}$ we will have then the equation  $\mu
\lbrack p_{0}- \frac{p_{1}^{2}}{4}]=1$ $\Rightarrow \mu
(p_{1}^{2}-4p_{0})=-4$, since the fact that $p_{1}^{2}-4p_{0}$ is
not a square $P$ is  irreducible in $A$, we will have -$\mu $ no
square, that means that the number of possible  $\mu $ is
$\frac{q-1}{2}$, if we consider the fact that $p_{0},p_{1}$ are
fixe, we will have $\frac{(q-1)}{2}$ solutions, then it remains
to calculate : for the case $c \neq 2$, for this, we put : $\
i_{2}=T+j_{2}$, $j_{2}\in \mathbf{F}_{q}$ and $c=aT+b$ where $a\in
\mathbf{F}_{q}^{\ast }$ and $b\in \mathbf{F}_{q}$.\ The fact that
$i_{2}\mid (c-2)\Rightarrow j_{2}=\frac{b-2}{a} $ and since
$i_{2}^{2}\mid P_{\Phi }(1)$ we will have : $1-(aT+b)+\mu
(T^{2}+p_{1}T+p_{0})=\mu (T+j_{2})^{2}\Rightarrow 1-b+\mu
p_{0}=\mu j_{2}^{2}$ \ \ and $\mu p_{1}-a=2\mu j_{2}$, then :
$\mu =\frac{a}{
p_{1}-2j_{2}}=\frac{a}{p_{1}-2(\frac{b-2}{a})}$ and then :%
\begin{equation*}
\frac{a}{p_{1}-2(\frac{b-2}{a})}[p_{0}-(\frac{b-2}{a})^{2}]+1-b=0\text{.}
\end{equation*}
The solution numbers of this equation in $(a,b)$ ( $p_{0},p_{1}$
are fixe ) give us

$\#\{i_{2},i_{2}^{2}\mid P_{\Phi }(1)$ et $i_{2}\mid (c-2)\},$
then we will have $(q-1)$ possible cases for $i_{2}$. Then :
\begin{eqnarray*}
N_{0}(2,1,q) &=&\frac{\#\{i_{2},i_{2}^{2}\mid P_{\Phi
}(1)\}}{\#\{\Phi ,
\text{ isogeny, ordinary}\}} \\
&=&\frac{(q-1)+\frac{(q-1)}{2}}{(q-1)[(\frac{q-1}{2})q-1]} \\
&=&\frac{\frac{3(q-1)}{2}}{(q-1)[(\frac{q-1}{2})q-1]} \\
&=&\frac{3}{q(q-1)-2}\Rightarrow
\end{eqnarray*}
\begin{eqnarray*}
C_{0}(2,1,q) &=&1-\frac{3}{q(q-1)-2} \\
&=&\frac{q(q-1)-5}{q(q-1)-2}.
\end{eqnarray*}
And for $N(2,1,q):$
\begin{equation*}
N(2,1,q)=\frac{\sum\limits_{P_{\Phi }}\sum\limits_{i_{2}}n(\Phi
,i_{2}). \text{\#}\{i_{2},i_{2}^{2}\mid P_{\Phi }(1)\text{ and
}i_{2}\mid (c-2)\}}{ \#\{\Phi \text{, isomorphism, ordinary}\}}
\end{equation*}
\begin{equation*}
=\frac{\sum\limits_{P_{\Phi
}}\sum\limits_{i_{2},i_{2}^{2}\shortmid 4-4\mu P}H(O(\frac{4-4\mu
P}{i_{2}^{2}}).\frac{q-1}{2}+\sum\limits_{P_{\Phi
}}\sum\limits_{i_{2},i_{2}^{2}\shortmid c^{2}-4\mu
P}H(O(\frac{c^{2}-4\mu P }{i_{2}^{2}}))(q-1)}{q^{3}-q^{2}-q+1}
\end{equation*}
\begin{equation*}
=\frac{\frac{q-1}{2}\sum\limits_{P_{\Phi
}}\sum\limits_{i_{2,}i_{2}^{2}\shortmid 4-4\mu P}H(O(\frac{4-4\mu
P}{ i_{2}^{2}}))+(q-1)\sum\limits_{P_{\Phi
}}\sum\limits_{i_{2},i_{2}^{2}\shortmid c^{2}-4\mu
P}H(O(\frac{c^{2}-4\mu P }{i_{2}^{2}}))}{q^{3}-q^{2}-q+1}.
\end{equation*}%
Finally : $C(2,1,q)=1-N(2,1,q)=$
\begin{equation*}
1-\frac{\frac{q-1}{2}\sum\limits_{P_{\Phi
}}\sum\limits_{i_{2,}i_{2}^{2}\shortmid 4-4\mu P}H(O(\frac{4-4\mu
P}{ i_{2}^{2}}))+(q-1)\sum\limits_{P_{\Phi
}}\sum\limits_{i_{2},i_{2}^{2}\shortmid c^{2}-4\mu
P}H(O(\frac{c^{2}-4\mu P }{i_{2}^{2}}))}{q^{3}-q^{2}-q+1}=
\end{equation*}

$\frac{q^{3}-q^{2}-q+1-[\frac{q-1}{2}\sum\limits_{P_{\Phi
}}\sum\limits_{i_{2,}i_{2}^{2}\shortmid 4-4\mu P}H(O(\frac{4-4\mu
P}{ i_{2}^{2}})+(q-1)\sum\limits_{P_{\Phi
}}\sum\limits_{i_{2},i_{2}^{2}\shortmid c^{2}-4\mu
P}H(O(\frac{c^{2}-4\mu P }{i_{2}^{2}}))]}{q^{3}-q^{2}-q+1}.$
\end{proof}

\section{the case : $m=2$ and $d=1$}

In this case $ n = m.d = 2 $, and $n_{2}=2\Rightarrow \#\{\Phi ,$
isomorphism, ordinary$\}=q^{3}-q-(q^{2}-1)=q^{3}-q^{2}-q+1$.

\begin{proposition}
Let $L=F_{q^{n}}$ and $P$ the $A$-characteristic of $L$,
$m=[L,A/P]$ and $\ d= $ deg $P$. We suppose that $m=2$ and $d=1$.
Then:

\begin{equation*}
C_{0}(1,2,q)=\frac{(q-1)q-4}{(q-1)q-2}\text{,}
\end{equation*}

\begin{equation*}
C(1,2,q)=\frac{q^{3}-q^{2}-q+1-\sum\limits_{P_{\Phi
}}\sum\limits_{i_{2,}i_{2}^{2}\shortmid c^{2}-4\mu
P}H(O(\frac{c^{2}-4\mu P }{i_{2}^{2}}))}{q^{3}-q^{2}-q+1}.
\end{equation*}
\end{proposition}
\begin{proof}
 We put $i_{2}=T+j_{2,}$ $j_{2}\in \mathbf{F}_{q}$ and $P(T)=T+p$ where
$p\in \mathbf{F}_{q}$. We start by calculate:
\begin{equation*}
N_{0}(1,2,q)=\frac{\#\{i_{2},i_{2}^{2}\mid P_{\Phi
}(1)\}}{\#\{\Phi ,\text{ isogeny, ordinary}\}},
\end{equation*}
In this case $c=2$, $i_{2}^{2}\mid P_{\Phi }(1)=\mu P^{2}-1$ we
will have  : $ 2\mu p=2\mu j_{2\text{ }}\mu p^{2}-1=\mu
j_{2}^{2}$ that means that $ p=j_{2}$ and $\mu
(p^{2}-j_{2}^{2})=1$ so the contradiction, and then $
\#\{i_{2},i_{2}^{2}\mid P_{\Phi }(1)\}=0$. For $c\neq 2$ we
calculate  $ \#\{i_{2},i_{2}^{2}\mid P_{\Phi }(1)$ and $i_{2}\mid
(c-2)\}$ and we remark first : $i_{2}\mid (c-2)\Rightarrow
j_{2}=\frac{b-2}{a}$ and $ i_{2}^{2}\mid P_{\Phi }(1)$ implies
that : $2p\mu -a=2\mu j_{2}$ and $ p^{2}+1-b=\mu j_{2}^{2}$,
finally we will have the  equation :
\begin{equation*}
p^{2}+1-b=\frac{a}{2p-2(\frac{b-1}{a})}(\frac{b-2}{a})^{2}
\end{equation*}
which is an equation in $(a,b,p)$ and accept  $(q-1)$ solutions,
since the fact that $p$ is  fix, then :%
\begin{eqnarray*}
N_{0}(1,2,q) &=&\frac{\#\{i_{2},i_{2}^{2}\mid P_{\Phi
}(1)\}}{\#\{\Phi \text{
; isogeny, ordinary}\}} \\
&=&\frac{q-1}{(q-1)((\frac{q-1}{2})q-1)}
\end{eqnarray*}%
\begin{eqnarray*}
&=&\frac{1}{(\frac{q-1}{2})q-1}\text{;} \\
&\Rightarrow &C_{0}(1,2,q)=1-N_{0}(1,2,q) \\
&=&1-\frac{1}{(\frac{q-1}{2})q-1)} \\
&=&\frac{(q-1)q-4}{(q-1)q-2}.
\end{eqnarray*}
And for $N(1,2,q)$ :
\begin{equation*}
N(1,2,q)=\frac{\sum\limits_{P_{\Phi }}\sum\limits_{i_{2}}n(\Phi
,i_{2}). \text{\#}\{i_{2},i_{2}^{2}\mid P_{\Phi }(1)\text{ and
}i_{2}\mid (c-2)\}}{ \#\{\Phi \text{, isomorphism, ordinary}\}}
\end{equation*}
\begin{equation*}
=\frac{\sum\limits_{P_{\Phi
}}\sum\limits_{i_{2,}i_{2}^{2}\shortmid c^{2}-4\mu
P}H(O(\frac{c^{2}-4\mu P}{i_{2}^{2}}))(q-1)}{q^{3}-q}
\end{equation*}
\begin{equation*}
=\frac{\sum\limits_{P_{\Phi
}}\sum\limits_{i_{2,}i_{2}^{2}\shortmid c^{2}-4\mu
P}H(O(\frac{c^{2}-4\mu P}{i_{2}^{2}}))}{q(q+1)}.
\end{equation*}
In end:%
\begin{eqnarray*}
C(1,2,q) &=&1-\frac{\sum\limits_{P_{\Phi
}}\sum\limits_{i_{2,}i_{2}^{2}\shortmid c^{2}-4\mu
P}H(O(\frac{c^{2}-4\mu P
}{i_{2}^{2}}))}{q^{3}-q^{2}-q-1}= \\
&&\frac{q^{3}-q^{2}-q+1-\sum\limits_{P_{\Phi
}}\sum\limits_{i_{2,}i_{2}^{2}\shortmid c^{2}-4\mu
P}H(O(\frac{c^{2}-4\mu P }{i_{2}^{2}}))}{q^{3}-q^{2}-q+1}\text{.}
\end{eqnarray*}
\end{proof}

\subsection{$\lim_{q\rightarrow \infty }C_{0}(d,m,q)$ and $\lim_{q\rightarrow
\infty }C(d,m,q)$ for $m.d$ $\leq 2$}

By the calculus of $C_{0}(d,m,q)$ and $C(d,m,q)$ for $m.d$ $\leq
2$, we have :

\begin{Corollary}
Let $L=F_{q^{n}}$ and $P$ the $A$-characteristic of $L$,
$m=[L,A/P]$ and $\ d=$ deg$P$. Then :
\end{Corollary}

\begin{equation*}
\lim\limits_{q\rightarrow \infty
}C_{0}(1,1,q)=\lim\limits_{q\rightarrow \infty
}C_{0}(1,2,q)=\lim\limits_{q\rightarrow \infty
}C_{0}(2,1,q)=1\text{,}
\end{equation*}

\begin{equation*}
\lim\limits_{q\rightarrow \infty
}C(1,1,q)=\lim\limits_{q\rightarrow \infty
}C(1,2,q)=\lim\limits_{q\rightarrow \infty }C(2,1,q)=1\text{.}
\end{equation*}

\begin{proof}
 Since : $C_{0}(1,2,q)=\frac{q(q-1)-4}{q(q-1)-2}$, $C_{0}(2,1,q)=
\frac{q(q-1)-5}{q(q-1)-2}$ et $C_{0}(1,1,q)=1$, where we can see
that : $md\leq 2$ , $\lim\limits_{q\rightarrow \infty
}C_{0}(d,m,q)=1$. By other way, since $C_{0}(d,m,q)$ $\leq $
$C(d,m,q)\leq 1$ and by passing to the limit, we have  :
$\lim\limits_{q\rightarrow \infty }C(d,m,q)=1$, pour $md\leq 2$.
\end{proof}

By the results above, we can give the following conjecture :

\begin{Conjecture}
Let $L=F_{q^{n}}$ and $P$ the $A$-characteristic of $L$,
$m=[L,A/P]$ and $\ d= $ deg$P$. Then :

\begin{equation*}
\lim\limits_{q\rightarrow \infty
}C(d,m,q)=\lim\limits_{q\rightarrow \infty }C_{0}(d,m,q)=1\text{.}
\end{equation*}
\end{Conjecture}
\subsection{Discussion and open Questions }

By the last conjecture we can ask whether for a big values of $d$
and $m$, we will have a cyclic modules ? otherwise  :

$\lim\limits_{d\rightarrow \infty
}C(d,m,q)=\lim\limits_{d\rightarrow \infty }C_{0}(d,m,q)=1$?

of same : $\lim\limits_{m\rightarrow \infty
}C(d,m,q)=\lim\limits_{m\rightarrow \infty }C_{0}(d,m,q)=1?$

We can also ask whether  the rank of a Drinfeld $A$-module $%
\Phi $ is  decisive for the cyclicity of  A-module $L^{\Phi }$?

Lastly, it is legitimate to ask if we can generalize the Theorem
4.1 for a Drinfeld $A$-modules $\Phi $ such that  $A$ is not
$F_{q}[T]$ and such that the rank of $\Phi $ is more bigger than
2?

\providecommand{\bysame}{\leavevmode\hbox to3em{\hrulefill}\thinspace} %
\providecommand{\MR}{\relax\ifhmode\unskip\space\fi MR } \providecommand{%
\MRhref}[2]{} \providecommand{\href}[2]{#2}

\end{document}